\documentclass[12pt]{article}
\usepackage{amsfonts}

\usepackage{latexsym}
\usepackage{amssymb}
\usepackage{amsmath}
\usepackage{amsthm}
\usepackage[mathscr]{eucal}
\usepackage{graphicx}
\usepackage{hyperref}
\usepackage{caption}
\usepackage{subcaption}

\newtheorem{Lemma}{Lemma}

\newtheorem{Theorem}{Theorem}

\newtheorem{Definition}{Definition}

\renewcommand{\qed}{\hfill{\ \ \rule{2mm}{2mm}} \vspace{0.2in}}

\begin{document}

\title{Critical Exponent for the Acyclic Chromatic Number of Random Graphs}
\author{ \textbf{Ghurumuruhan Ganesan}
\thanks{E-Mail: \texttt{gganesan82@gmail.com} } \\
\ \\
IISER Bhopal, India}
\date{}
\maketitle


\begin{abstract}
In this paper we study acyclic colouring in the random subgraph~\(G\) of the complete graph~\(K_n\) on~\(n\) vertices where each edge is present with probability~\(p,\) independent of the other edges. We show that the acyclic chromatic number exhibits a phase transition from sub-linear to linear growth as the edge probability increases, even in the sparse regime and obtain estimates for the critical exponent.  Next, we introduce a relaxation by allowing for a small fraction of ``bad" cycles to violate the acyclic colouring condition and show that the critical exponent in this case is in fact zero, no matter how small the fraction. 

\vspace{0.1in} \noindent \textbf{Key words:} Acyclic colouring; random graphs; critical exponent; fractional acyclic chromatic number.

\vspace{0.1in} \noindent \textbf{AMS 2000 Subject Classification:} Primary: 05C62;
\end{abstract}

\bigskip

\renewcommand{\theequation}{\thesection.\arabic{equation}}
\setcounter{equation}{0}
\section{Introduction} \label{intro}

Acyclic chromatic number~\(\chi_{acyc}(H)\) of a deterministic graph~\(H\) is the smallest number of colours needed for obtaining a proper colouring of a graph with no bichromatic cycles. Such colourings are important from both theoretical and application perspectives and there are well studied bounds on the acyclic chromatic number. Alon et al. (1991) use the Local Lemma to see that~\(\chi_{acyc}(H)\)  grows at most of the order of~\(\Delta^{\frac{4}{3}}\) where~\(\Delta\) is the maximum vertex degree of~\(G\) and also provide examples of graphs that nearly achieve the bound. Recently, Sereni and Volec use entropy compression techniques (Esperet and Parreau (2013)) to obtain improvements on the constants in the bounds derived in Alon et al. (1991).  In a separate line of work, there are studies that determine \emph{exactly} the acyclic chromatic number for graphs with small vertex degrees (see Yadav et al. (2011) and references therein). For more on the application of probabilistic techniques in graph colourings we refer to the book by Molloy and Reed (2002).

In this paper, we study acyclic colouring in random graphs. We obtain estimates on the critical exponent at which the acyclic chromatic number transitions from sublinear to linear growth (with respect to the number of vertices~\(n\)). We then show that by allowing for an arbitrary small fraction of cycles to violate the acyclic condition, the resulting ``fractional" acyclic chromatic number always remains sublinear in the number of vertices~\(n.\)







The paper is organized as follows. In Section~\ref{sec_app_two} we state and prove our main result regarding estimates for the critical exponent of the acyclic chromatic number in random graphs.

\setcounter{equation}{0}
\renewcommand\theequation{\arabic{section}.\arabic{equation}}
\section{Acyclic Chromatic Number}\label{sec_app_two}
Let~\(K_n\) be the complete graph on~\(n\) vertices and let~\(X_f, f \in K_n\) be independent Bernoulli random
variables indexed by the edge set of~\(K_n,\) with distribution
\begin{equation}\label{x_dist}
\mathbb{P}(X_f = 1)= p(f)  = 1-\mathbb{P}(X_f = 0)
\end{equation}
where~\(0 \leq p(f) \leq 1.\) Let~\(G\) be the random subgraph of~\(K_n\) formed by the set of all edges satisfying~\(X_f = 1.\) If~\(p(f)= p\) for all edges~\(f,\) then~\(G\) is said to be \emph{homogenous}; otherwise we say that~\(G\) is \emph{inhomogenous}.

A path in~\(G\) is a sequence of vertices~\(P = (u_1,\ldots,u_t)\) such that~\(u_i\) is adjacent to~\(u_{i+1}\) for each~\(1 \leq i \leq t-1.\) The length of~\(P\) is the number of edges in~\(P\) which in this case is~\(t-1.\) If in addition, the vertex~\(u_t\) is also adjacent to~\(u_1\) in the graph~\(G,\) then~\(P\) is said to be a \emph{cycle} of length~\(t.\)


For an integer~\(k \geq 1,\) a proper~\(k-\)colouring of~\(G\) is a map~\(g: V \rightarrow \{1,2,\ldots,k\}\) such that~\(g(u) \neq g(v)\) for all edges~\((u,v) \in E.\) The chromatic number~\(\chi(G)\) is the smallest integer~\(k\) such that~\(G\) admits a proper~\(k-\)colouring.  We have the following definition.
\begin{Definition}\label{def_one}
Let~\(f\) be a proper colouring of~\(G.\) For~\(0 < \eta < 1\) and integer~\(c \geq 3,\) we say that~\(f\) is~\((\eta,c)-\)acyclic if a fraction of at least~\((1-\eta)\) of the all the cycles of length~\(c\) or more in~\(G,\) contain at least~\(c\) distinct colours.

We define~\(\chi_{\eta,c}(G),\) the minimum number of colours needed for a proper\\\((\eta,c)-\)acyclic colouring of~\(G,\) to be the~\((\eta,c)-\)\emph{acyclic chromatic number} of~\(G.\)
\end{Definition}
For~\(\eta=0,c=3\) and~\(p(f) = p\) for all edges~\(f \in K_n,\) the above definition reduces to the usual acyclic chromatic number as in Alon et al. (1991). In our definition above, we allow for a small fraction~\(\eta\) of the cycles to possibly violate the acyclic colouring condition and study the effect of this relaxation in Theorem~\ref{thm_main} below.  



If~\(p(f) = p\) for all edges~\(f,\) then both the chromatic number and the acyclic chromatic number increase with the edge probability~\(p.\) However, unlike the chromatic number~\(\chi(G)\) that remains sub-linear even when the edge probability~\(p\) becomes a constant (Alon and Spencer (2008)), the acyclic chromatic number exhibits a ``transition" from sublinear to linear growth even in the sparse regime, where~\(p\) is of the form~\(\frac{1}{n^{\beta}}\) for some constant~\(\beta > 0.\) We illustrate this aspect below by considering the case~\(c=3\) as an example. Throughout, we use the following standard deviation estimate regarding sums of Bernoulli random variables. Let~\(Z_i, 1 \leq i \leq t\) be independent Bernoulli random variables satisfying~\(\mathbb{P}(Z_i = 1) = p_i = 1-\mathbb{P}(Z_i = 0).\) If~\(W_t = \sum_{i=1}^{t} Z_i\) and~\(\mu_t = \mathbb{E}W_t,\) then for any~\(0 < \epsilon < \frac{1}{2}\) we have that
\begin{equation}\label{conc_est_f}
\mathbb{P}\left(\left|W_t-\mu_t\right| \geq \epsilon \mu_t\right) \leq 2\exp\left(-\frac{\epsilon^2}{4}\mu_t\right).
\end{equation}
For a proof of~(\ref{conc_est_f}), we refer to Corollary~\(A.1.14,\) pp.~\(312,\) Alon and Spencer (2008).

We recall from Theorem~\(1.1\) of Alon et al. (1991) that for any deterministic graph~\(\Gamma\) the acyclic chromatic number~\(\chi_{0,3}(\Gamma) \leq C\Delta_{\Gamma}^{4/3},\) where~\(\Delta_{\Gamma} \) is the maximum vertex degree of~\(\Gamma\) and~\(C > 0\) is an absolute constant not depending on~\(\Gamma\) or~\(n.\) To translate this bound to~\(G,\) we note that the expected degree of any vertex in~\(G\) equals~\((n-1)p \leq np.\) Therefore from~(\ref{conc_est_f}) and the union bound, we get that if~\(p = \frac{1}{n^{\beta}}\) then~\(\Delta \leq 2np = 2n^{1-\beta}\) with high probability, i.e. with probability converging to one as~\(n \rightarrow \infty.\) Consequently,~\(\frac{\chi_{0,3}(G)}{n} \leq \frac{C}{n}\Delta^{4/3} \leq \frac{C}{n} (2np)^{4/3}  = C2^{4/3} n^{\frac{(1-4\beta)}{3}} \longrightarrow 0\) in probability, if we set~\(\beta > \frac{1}{4}\) strictly.  On the other hand from the proof of Theorem~\(1.2\) of Alon et al. (1991),  we also know that if~\(\beta < \frac{1}{4}\) then~\(\mathbb{P}\left(\chi_{0,3}(G) \geq \frac{n}{4}\right) \longrightarrow 1\) as~\(n \rightarrow \infty.\)  Thus~\(\beta = \frac{1}{4}\) could be interpreted as the ``critical exponent" at which the acyclic chromatic number transitions from being sub-linear to growing linearly with~\(n.\)



Generalizing the above discussion, let~\(H = H_n \subset K_n\) be any sequence of deterministic graphs and  let~\(G\) be the random graph as in~(\ref{x_dist}) with edge probability~\(p(f) = \frac{1}{n^{\beta}} =: p\) for all~\(f \in K_n \setminus H\) and~\(p(f)= 0\) otherwise. Define
\begin{equation}\label{beta_crit_def}
\beta_{crit}(\eta, c) := \inf\left\{\beta > 0 : \frac{\chi_{\eta,c}(G)}{n}\longrightarrow 0 \text{ in probability}\right\}
\end{equation}
to be the \emph{critical exponent} for the~\((\eta,c)-\)acyclic chromatic number. That\\\(\beta_{crit}(\eta,c)\) is finite is obtained from the following direct edge counting argument: if~\(m(G)\) is the number of edges in~\(G,\) then~\(\mathbb{E}m(G) \leq n^2p\) and so by the Markov inequality, we know that~\(\mathbb{P}(m(G) \geq n^2p\log{n}) \leq \frac{1}{\log{n}}.\) Setting~\(p = \frac{1}{n^3},\) we get that~\(m(G) =0\) and therefore~\(\chi_{\eta, c}(G) = 0\) with high probability. This implies that~\(\beta_{crit}(\eta,c) < 3.\)



The following result obtains sharper bounds for the critical exponent in terms of~\(c.\) Throughout constants do not depend on~\(n\) and for a constant~\(0 < s_0 < 1,\) we say that a deterministic graph~\(H \subset K_n\) has edge density at most~\(s_0\)  if for any~\(1 \leq l \leq n\) and any set~\({\cal V}\) containing~\(l\) vertices, the number of edges of~\(H\) with both endvertices in~\({\cal V}\) is at most~\(s_0 {l \choose 2}.\)
\begin{Theorem}\label{thm_main} Suppose~\(H = H_n\) is a sequence of deterministic graphs each with edge density strictly less than~\(\frac{1}{4}.\)\\
\((a)\) For every integer constant~\(c \geq 3,\) we have that~\[ \frac{1}{4c-2} \leq \beta_{crit}(0,c) \leq 1-\frac{1}{c-1}.\]
\((b)\) For every~\(0 < \eta < 1\) and~\(c \geq 3,\) we have that~\(\beta_{crit}(\eta,c) = 0.\)
\end{Theorem}
From part~\((a),\) we see that for the  ``strict" acyclic colouring with~\(\eta = 0,\) the critical exponent is non-trivial, even if~\(H\) has positive edge density; i.e. even if a fraction of the total edges of~\(K_n\) are removed. For~\(\eta > 0\) however, part~\((b)\) implies that acyclic chromatic number is \emph{always} sublinear in~\(n\) with high probability, no matter how small~\(\eta, \beta\) or~\(H\) is. In other words, by allowing a small fraction of cycles to violate the acyclic condition, we obtain the desired colouring using~\(o(n)\) colours with high probability.





To prove Theorem~\ref{thm_main}\((a),\) we use the following auxiliary result regarding paths of given length in~\(G.\) Letting~\(F_k\) be the event that there exists a path of length~\(k\) in~\(G,\) we have the following result.
\begin{Lemma}\label{lem_path} Let~\(H \subset K_n\) be any deterministic graph and suppose~\(p(f) = p\) for all edges~\(f \in K_n \setminus H\) and~\(p(f) = 0\) otherwise. If~\(H\) has edge density at most~\(s_0 < 1,\) then
\begin{equation}\label{cyc_est_ax22}
\mathbb{P}(F_k) \geq 1-k n^{k-1} e^{-Cn^2p^{2k-1}}
\end{equation}
for some constant~\(C  > 0.\) 
\end{Lemma}

\emph{Proof of Lemma~\ref{lem_path}}: Defining~\(t_j(l) \) to be the maximum probability that any induced subgraph~\(G\) containing~\(l\) vertices \emph{does not} contain a path of length~\(j,\) we obtain a recursion involving~\(t_k(n).\) First, if~\(m(G)\) is the number of edges in~\(G,\) we know that~\( {n \choose 2}p \geq \mathbb{E}m(G) \geq  {n \choose 2}(1-s_0)p\) and so from the concentration estimate~(\ref{conc_est_f}), we get that
\begin{equation}\label{dimple2}
\mathbb{P}\left(m(G) \geq \frac{n^2p(1-s_0)}{4}\right) \geq 1-e^{-Dn^2p}
\end{equation}
for some constant~\(D > 0.\)

Next, if~\(m(G) \geq \frac{n^2p(1-s_0)}{4},\) then using the fact that the sum of vertex degrees is twice the number of edges, there must exist a vertex~\(z\) whose degree~\(d(z)\) in~\(G\) is at least~\(\frac{np(1-s_0)}{2}.\) Thus~\(t_k(n) = \mathbb{P}(F_k^c)\) is bounded above as
\begin{eqnarray}
t_k(n) &\leq& \mathbb{P}\left(F_k^c \bigcap \bigcup_{1 \leq z \leq n} \left\{d(z) \geq \frac{np(1-s_0)}{2}\right\}\right)  + \mathbb{P}\left(m(G) \leq \frac{n^2p(1-s_0)}{4}\right) \nonumber\\
&\leq& \sum_{1 \leq z \leq n}\mathbb{P}\left(F_k^c \bigcap \left\{d(z) \geq \frac{np(1-s_0)}{2}\right\}\right)+e^{-Dn^2p} \label{dimple}
\end{eqnarray}
by~(\ref{dimple2}).

To estimate the first term in the last expression of~(\ref{dimple}), we let~\({\cal N}(z)\) be the set of neighbours of~\(z\) in~\(G\)  and write
\begin{equation}\label{dimple3}
\mathbb{P}\left(F_k^c \bigcap \left\{d(z) \geq \frac{np(1-s_0)}{2}\right\}\right) = \sum_{S} \mathbb{P}\left(F_k^c \cap \{{\cal N}(z) = S\}\right),
\end{equation}
where the summation is over all sets~\(S\) of cardinality at least~\(\frac{np(1-s_0)}{2}\) and not containing the vertex~\(z.\) If~\(F_k^c \cap \{{\cal N}(z) = S\}\) occurs then so does~\(R_{k-1}(S) \cap \{{\cal N}(z) = S\},\) where~\(R_{k-1}(S)\) is the event that there is no path of length~\(k-1\) in the induced subgraph of~\(G\) with vertex set~\(S.\) The event~\(R_{k-1}(S)\) is independent of~\(\{{\cal N}(z) = S\}\) and so
\begin{eqnarray}
\mathbb{P}\left(F_k^c \cap \{{\cal N}(z) = S\}\right) &\leq& \mathbb{P}\left(R_{k-1}(S) \cap \{{\cal N}(z) = S\}\right)  \nonumber\\
&=& \mathbb{P}(R_{k-1}(S)) \mathbb{P}({\cal N}(z) = S) \nonumber\\
&\leq& t_{k-1}(n_1) \mathbb{P}({\cal N}(z) = S), \nonumber
\end{eqnarray}
since~\(S\) has size at least~\(n_1 := \frac{np(1-s_0)}{2}.\) Plugging this into~(\ref{dimple3}) we get
\[\mathbb{P}\left(F_k^c \bigcap \left\{d(z) \geq \frac{np(1-s_0)}{2}\right\}\right) \leq t_{k-1}(n_1) \sum_{S} \mathbb{P}({\cal N}(z) = S) \leq t_{k-1}(n_1)\] and so from~(\ref{dimple}), we get that
\begin{equation}\label{imp_rec}
t_k(n) \leq nt_{k-1}(n_1) + e^{-Dn^2p}.
\end{equation}

Setting~\(n_i := n \left(\frac{p}{2}(1-s_0)\right)^{i}\) for~\(i \geq 0\) and applying~(\ref{imp_rec}) recursively, we get that
\begin{eqnarray}
t_k(n) &\leq& n(n_1 t_{k-2}(n_2) + e^{-Dn_1^2p}) + e^{-Dn^2p}  \nonumber\\
&\leq& n^2t_{k-2}(n_2) + ne^{-Dn_1^2p} + e^{-Dn^2p} \nonumber\\
&\leq& \ldots \nonumber\\
&\leq& n^{k-1} t_1(n_{k-1}) + \sum_{i=0}^{k-2}n^{i}e^{-Dn_i^2p}. \nonumber
\end{eqnarray}
From~(\ref{dimple2}), we see that~\(t_1(n_{k-1}) \leq e^{-Dn_{k-1}^2p}\) and so
\[t_k(n) \leq \sum_{i=0}^{k-1}n^{i}e^{-Dn_i^2p} \leq k n^{k-1}e^{-Dn_{k-1}^2p} \leq k n^{k-1} e^{-Cn^2p^{2k-1}}\]
for some constant~\(C > 0.\) This completes the proof of Lemma~\ref{lem_path}.~\(\qed\)

To prove Theorem~\ref{thm_main}\((b)\) we use the following version of the Lov\'asz Local Lemma (see Lemma~\(2.1\) of Alon et al. (1991)).
\begin{Lemma}\label{local} Let~\(A_1,\ldots,A_t\) be events in an arbitrary probability space. Let~\(\Gamma \) be the dependency graph for the events~\(\{A_i\},\) with vertex set~\(\{1,2,\ldots,t\}\) and edge set~\({\cal E};\) i.e. assume that each~\(A_i\) is independent of the family of events~\(A_j, (i,j) \notin {\cal E}.\) If there are reals~\(0 \leq y(i) < 1\) such that~\(\mathbb{P}(A_i) \leq y(i) \prod_{(i,j)\in {\cal E}} (1-y(j)),\) for each~\(i,\) then~\[\mathbb{P}\left(\bigcap_{i} A^c_i\right) \geq \prod_{1 \leq i \leq n} (1-y(i)) > 0.\]
\end{Lemma}

We now use Lemma~\ref{lem_path} to prove Theorem~\ref{thm_main}\((a)\) below.\\
\emph{Proof of Theorem~\ref{thm_main}\((a)\)}: For the upper bound, we assume that~\(H = \emptyset\) and estimate~\(\chi_{0,c}(G)\) using the power graph of~\(G\) as follows. The~\(k^{th}\) power of~\(G\) is the graph~\(G_k\) obtained by connecting two vertices~\(u\) and~\(v\) by an edge if the distance between~\(u\) and~\(v\) in~\(G\) is at most~\(k.\) Any proper colouring of~\(G_{c-1}\) is a proper~\((0,c)-\)acyclic colouring of~\(G\) and so~\(\chi_{0,c}(G) \leq \Delta(G_{c-1})+1\) where~\(\Delta(G_{c-1})\) is the maximum vertex degree in~\(G_{c-1}.\) By construction, if~\(\Delta_G\) is the maximum vertex degree of~\(G,\) then~\(\Delta(G_{c-1}) \leq \Delta_G^{c-1}\) and since~\(np  = n^{1-\beta} \geq (\log{n})^2,\) we get from the deviation estimate~(\ref{conc_est_f}) and the discussion in the paragraph prior to~(\ref{beta_crit_def}) that~\(\chi_{c}(G) \leq \Delta_G^{c-1} \leq (2np)^{c-1}\) with high probability. Choosing~\(\beta  > 1-\frac{1}{c-1}\) strictly, we then get that~\(\frac{\chi_{0,c}(G)}{n} \longrightarrow 0\) in probability and so~\(\beta_{crit}(0,c) \leq 1-\frac{1}{c-1}.\)


For the lower bound on~\(\beta_{crit}(0,c),\) we set~\(p = \frac{1}{n^{\beta}}\) with~\(\beta < \frac{1}{4c-2}\) strictly and show that with high probability, no~\(t-\)proper colouring of~\(G\) for~\(t \leq \frac{n}{4},\) contains a cycle of length~\(c\) with at most~\(c-1\) distinct colours. For simplicity we assume that~\(c\) is even henceforth (an analogous analysis holds otherwise with~\(\frac{c}{2}\) replaced by~\(\frac{c-1}{2}\) below) and  begin with some preliminary computations. Let~\({\cal V}=\{V_1,\ldots,V_t\}\) be any partition of~\(\{1,2,\ldots,n\}\) and say that~\({\cal V}\) is a \emph{good} partition if~\({\cal V}\) contains~\(z \geq \frac{np}{100\log{n}}\) sets, each containing at least two vertices. Suppose~\({\cal V}\) is a good partition and let~\(U_1,\ldots,U_z\) be the sets in~\({\cal V}\) each containing at least two vertices. Pick two vertices~\(u_{i,1}\) and~\(u_{i,2}\) each from each~\(U_i\) according to a deterministic rule and consider the graph~\(K({\cal V})\) with vertex set~\(\{q_1,q_2,\ldots,q_z\}\) obtained as follows. The vertices~\(q_a\) and~\(q_b\) are joined by an edge if and only if~\(u_{a,j}\) is adjacent to both~\(u_{b,1}\) and~\(u_{b,2}\) in~\(K_n \setminus H\) for each~\(j=1,2.\)

Say that an edge~\((q_a,q_b)\) of~\(K({\cal V})\) is \emph{open} if each of the four edges in
\[{\cal E}_{a,b} := \{(u_{a,1},u_{a,2}),(u_{a,1},u_{b,1}), (u_{a,2},u_{b,1}), (u_{a,2},u_{b,2})\}\] is present in~\(G\) and let~\(G_{\cal V} \subset K({\cal V})\) be the set of all open edges. Letting~\(E_{path}\) be the event that~\(G_{\cal V}\) contains a path of length~\(\frac{c}{2}\) for \emph{every} good partition~\({\cal V},\) we compute the probability of the occurrence of~\(E_{path}\) as follows.  For a fixed good partition~\({\cal V},\) each edge of~\(K({\cal V})\) is open with probability~\(p^{4}\) and the total number of edges in the set~\(\bigcup_{1 \leq a < b \leq z} {\cal E}_{a,b}\) is~\(4 {z \choose 2}\) of which at most~\(4s_0 {z \choose 2}\) belong to~\(H.\) Thus there are at least~\((1-4s_0){z \choose 2}\) edges in~\(K({\cal V}).\)

We therefore argue as in Lemma~\ref{lem_path} with~\(n\) replaced by~\(z,\) the term~\(p\) replaced by~\(p^{4}\) and~\(s_0\) replaced by~\(4s_0,\) to get that~\(G_{\cal V}\) contains an open path of length~\(k=\frac{c}{2}\) with probability at least
\[1-z^{c/2}e^{-C_1 z^2p^{4(c-1)}} \geq 1-n^{c/2}\exp\left(-C_2\frac{n^2p^{4c-2}}{(\log{n})^2}\right) \geq 1-e^{-2Dn(\log{n})^2}\]
for some constants~\(C_1,C_2, D > 0\) and all~\(n\) large, provided~\(p = \frac{1}{n^{\beta}}\) with~\(\beta < \frac{1}{4c-2}\) strictly. Fixing such a~\(\beta\) and using the fact that there are at most~\(n^{n} = e^{n\log{n}}\) choices for~\({\cal V},\) we get from the union bound that
\begin{equation}\label{e_path_est}
\mathbb{P}(E_{path}) \geq 1-e^{n\log{n}}e^{-2Dn(\log{n})^2} \geq 1-e^{-Dn(\log{n})^2}.
\end{equation}


Our next step involves using~(\ref{e_path_est}) to arrive at a contradiction. Let~\({\cal V} = \{V_1,\ldots,V_t\}\) be a partition of the vertex set~\(\{1,2,\ldots,n\}\) representing a~\(t-\)proper colouring of~\(G\) so that the vertices in~\(V_i\) have colour~\(i.\) We denote~\(V_i\) to be a colour \emph{class} of~\({\cal V}\) and letting~\(E_{size}\) be the event that each colour class in a proper colouring of~\(G\) contains at most~\(\frac{8\log{n}}{p}\) vertices, we first show that~\(E_{size}\) occurs with high probability. Indeed, by definition, each~\(V_i\) must be a stable set and since the edge density of~\(H\) is at most~\(s_0,\) the probability that some~\(V_i, 1 \leq i \leq t\) has size~\(t \geq \frac{8\log{n}}{p}\) is
\begin{eqnarray}
\mathbb{P}(E_{size}^c) &\leq& {n \choose t} \cdot (1-p)^{{t \choose 2}(1-s_0)} \nonumber\\
&\leq& e^{t\log{n}}\exp\left(-\frac{p(1-s_0)t^2}{4}\right) \nonumber\\
&\leq& \exp\left(-\frac{p(1-s_0)t^2}{8}\right)  \label{aditi}
\end{eqnarray}
where the first inequality in~(\ref{aditi}) is true since~\({n \choose t} \leq n^{t} = e^{t\log{n}}\) and\\\( (1-p)^{t\choose 2} \leq e^{-p{t \choose 2}} \leq \exp\left(-\frac{pt^2}{4}\right)\) and the second inequality in~(\ref{aditi}) follows from the fact that~\(t \geq \frac{8\log{n}}{p}.\)


Suppose~\(E_{size}\) occurs so that each~\(V_i, 1 \leq i\leq t\) has size at most~\(\frac{8\log{n}}{p}.\) This necessarily implies that~\(t \geq \frac{np}{9\log{n}}\) and we now count the number of ``bad" sets in~\({\cal V}\) that have exactly one vertex in them. If~\(V_1,\ldots,V_x\) are the single vertex sets in~\({\cal V},\) then the remaining~\(n-x\) vertices must be ``packed" in the remaining~\(t-x\) sets~\(V_{x+1},\ldots,V_t.\) This must imply that at least one set in~\({\cal V}\) has size at least~\(\frac{n-x}{t-x} \geq \frac{10\log{n}}{p}\) if~\[x \geq \frac{t- \frac{np}{10\log{n}}}{1-\frac{p}{10\log{n}}} =: x_0\] and this leads to a contradiction. Thus~\(x \leq x_0\) and there are at least
\begin{equation}\label{y_not_est}
y_0 := t-x_0 \geq (n-t)\frac{p}{10\log{n}} \left(1-\frac{p}{10 \log{n}}\right)^{-1}
\end{equation}
good sets in~\({\cal V}\) each containing at least two vertices.

Summarizing the discussion in the previous two paragraphs we have the following: If~\(E_{bad}\) is the event that there exists a proper~\((0,c)-\)acyclic colouring of~\(G\) containing~\( t\leq \frac{n}{8}\) colours, then the occurrence of~\(E_{bad} \cap E_{size} \)  necessarily implies that there exists a proper~\((0,c)-\)acyclic colouring~\({\cal V} = \{V_1,\ldots,V_t\}\) containing~\(\frac{8\log{n}}{p} \leq t \leq \frac{n}{8}\) colour classes and  moreover, from~(\ref{y_not_est}), we also get that there at least
\begin{equation}\label{z_est}
\frac{7np}{80\log{n}} \left(1-\frac{p}{10 \log{n}}\right)^{-1} \geq \frac{7np}{100\log{n}} =: z
\end{equation}
colour classes, each containing at least two vertices.

\begin{figure}[tbp]
\centering
\includegraphics[width=3in, trim= 40 370 40 190, clip=true]{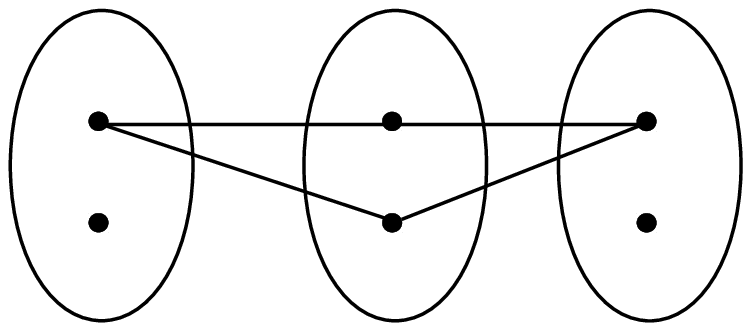}
\caption{Obtaining a cycle of length~\(c\) containing~\(\frac{c}{2}+1\) distinct colours for~\(c=4.\) The solid lines in the figure correspond to the cycle.}
\label{fig_1}
\end{figure}

Suppose now that~\(E_{bad} \cap E_{size} \cap E_{path}\) occurs where~\(E_{path}\) is the event defined prior to~(\ref{e_path_est}) and let~\({\cal V}\) be the proper~\((0,c)-\)acyclic colouring of~\(G\) using~\(t \leq \frac{n}{8}\) colours as described in the previous paragraph. From the description prior to~(\ref{e_path_est}), we then get that the graph~\(G_{\cal V}\) corresponding to~\({\cal V}\) contains a path of length~\(\frac{c}{2}\) and therefore~\(G\) contains a cycle of length~\(c\) with~\(\frac{c}{2}+1\) distinct colours. This is illustrated in Figure~\ref{fig_1} below for the case~\(c=4.\) This leads to a contradiction and so \[\mathbb{P}(E_{bad}) \leq \mathbb{P}(E_{path}^c \cup E_{size}^c) \leq \mathbb{P}(E_{path}^c) + \mathbb{P}(E_{size}^c) \longrightarrow 0\] as~\(n \rightarrow \infty\) by~(\ref{e_path_est}) and~(\ref{aditi}). This completes the proof of Theorem~\ref{thm_main}\((a).\)~\(\qed\)


Using Lemma~\ref{local}, we now prove Theorem~\ref{thm_main}\((b).\)\\
\emph{Proof of Theorem~\ref{thm_main}\((b)\)}: We set~\(H = \emptyset\) throughout and first obtain an upper bound for~\(\chi_{\eta,c}\) in two steps. In the first step, we show that for any integer~\(g \geq c+1,\) the number of cycles of length at most~\(g-1,\) is at most a fraction~\(\eta\) of the total number of cycles in~\(G.\) Next, we use Lemma~\ref{local} to estimate the number of colours needed to obtain a proper colouring of~\(G\) in which each path of length~\(g\) or more, contains at least~\(c\) distinct colours. Finally, we choose~\(g\) large enough so that~\(\chi_{\eta,c} = o(n)\) with high probability. We provide the details below.

Any cycle~\({\cal C}\) of length~\(k\) in~\(K_n\) is present in the graph~\(G\) with probability~\(p^{k}\) and there are at most~\(n^{k}\) and at least~\(\frac{1}{2k}n(n-1)\cdots(n-k+1) \geq \frac{n^{k}}{4k}\) choices for the vertex set of~\({\cal C}.\)  Thus if~\(N_k\) is the number of cycles of length~\(k\) in the random graph~\(G,\) then~\(\frac{(np)^{k}}{4k} \leq \mathbb{E}N_k \leq (np)^{k}\) and from Lemma~\(3.5\) of Janson et al. (2000), we also know that
\[var(N_k) \leq C_1 \frac{(\mathbb{E}N_k)^2}{n^2p} = o(\mathbb{E}N_k)^2\] for some constant~\(C_1 > 0.\) Defining~\(E_k := \left\{\frac{(np)^{k}}{8k} \leq N_k \leq 2(np)^{k}\right\},\) we get from the Chebychev's inequality that
\[\mathbb{P}\left(E_k\right) \geq 1 - \frac{C_2}{n^2p}\] for some constant~\(C_2 > 0.\) Letting~\(g \geq c+1\) and~\(E_{cyc} := \bigcap_{c \leq k \leq g+1}E_k,\) we get from the union bound that
\begin{equation}\label{e_cyc_est}
\mathbb{P}(E_{cyc}) \geq 1- \frac{C_3}{n^2p} \longrightarrow 1
\end{equation}
as~\(n \rightarrow \infty,\) for some constant~\(C_3 > 0.\) If~\(E_{cyc}\) occurs, then the total number of cycles of length at most~\(g-1\) is at most~\(\sum_{k=3}^{g-1} 2(np)^{k} \leq C_4 (np)^{g-1}\) for some constant~\(C_4 > 0\) and the number of cycles of length~\(g\) is at least~\(\frac{(np)^{g}}{4g}.\) Thus the number of cycles of length at most~\(g-1,\) is at most a fraction~\(\frac{4gC_4}{np} \leq \eta\) of the total number of cycles in~\(G.\) This completes the first step in our proof.

In the second step, we use Lemma~\ref{local} to obtain a proper colouring in which each \emph{path} of length~\(g\) or more has at least~\(c\) colours. This in turn ensures that each cycle of length~\(g\) or more has at least~\(c\) colours. Let~\(E_{deg}\) be the event that the degree of each vertex in~\(G\) lies between~\(\frac{np}{2}\) and~\(2np.\) Using the standard deviation estimate~(\ref{conc_est_f}) we have that
\begin{equation}\label{e_deg_est}
\mathbb{P}(E_{deg}) \geq 1- ne^{-Cnp}
\end{equation}
for some constant~\(C > 0.\) We henceforth assume that~\(E_{deg} \cap E_{cyc}\) occurs. For integer~\(\theta \geq 1\) to be determined later, let~\(X_1,\ldots,X_n\) be independent and uniformly distributed in~\(\{1,2,\ldots,\theta\},\) also independent of~\(G\) and let~\(\mathbb{P}_X\) denote the distribution of~\((X_1,\ldots,X_n).\) For~\(1 \leq u < v \leq n\) let~\(A_{u,v}\) be the event that~\(X_u = X_v.\)  For distinct integers~\({\cal U} := \{u_1,\ldots,u_g\},\) let~\(B_{{\cal U}}\) be the event that the set~\(\{X_{u_1}, \ldots, X_{u_g}\}\) has cardinality at most~\(c-1.\) There are constants~\(D_l, 1 \leq l \leq c-1\) such that
\begin{equation}\label{weight_ass}
\mathbb{P}_X(A_{u,v}) = \frac{1}{\theta} =: \frac{y(u,v)}{2} \text{ and } \mathbb{P}_X(B_{{\cal U}})\leq \sum_{l=1}^{c-1}{\theta \choose l}\frac{D_l}{\theta^{g}} \leq \frac{D \theta^{c-1}}{\theta^{g}} =: \frac{y({\cal U})}{2}.
\end{equation}

We now use Lemma~\ref{local} with weights~\(y(.)\) as defined above. We construct the dependency graph~\(\Gamma\) whose vertices are either edges of~\(G\) which we call as type~\(A\) vertex or paths of length~\(g\) in the graph~\(G\) which we call as type~\(B\) vertex. We recall that the maximum degree of any vertex is at most~\(2np.\) Therefore for any vertex~\(u \in K_n,\) there are at most~\(g (2np)^{g}\) paths in~\(G\) containing~\(u\) as an endvertex. This implies any type~\(A\) vertex in~\(\Gamma\) is adjacent to at most~\(2np\) type~\(A\) vertices and at most~\(2g (2np)^{g}\) type~\(B\) vertices. Similarly, any type~\(B\) vertex is adjacent to at most~\(2gnp\) type~\(A\) vertices and at most~\(g^2 (2np)^{g}\) type~\(B\) vertices. Consequently, if we ensure that
\begin{equation}\label{lova_a}
\frac{1}{\theta} \leq \frac{2}{\theta} \left(1-\frac{2}{\theta}\right)^{2np} \left(1-\frac{D\theta^{c-1}}{\theta^{g}}\right)^{2g(2np)^{g}}
\end{equation}
and
\begin{equation}\label{lova_b}
\frac{D \theta^{c-1}}{\theta^{g}}  \leq \frac{2D \theta^{c-1}}{\theta^{g}} \left(1-\frac{2}{\theta}\right)^{2gnp}\left(1-\frac{D\theta^{c-1}}{\theta^{g}}\right)^{g^2(2np)^{g}},
\end{equation}
then Lemma~\ref{local} with weight assignments as in~(\ref{weight_ass}) guarantees the existence of a proper colouring of~\(G\) using~\(\theta\) colours, in which every path of length~\(g\) contains at least~\(c\) distinct colours.

Setting~\(\theta = 8D \left(g^2 (2np)^{g}\right)^\frac{1}{g-c+1}\) and using~\((1-x)^{a} \geq 1-ax,\) we see that both~(\ref{lova_a}) and~(\ref{lova_b}) are satisfied for all~\(n\) large and so summarizing, we get that if~\(E_{deg} \cap E_{cyc}\) occurs, then~\(\chi_{\eta,c}(G) \leq  C (np)^{\frac{g}{g-c+1}}.\) Finally, choosing~\(g = g(\beta,c)\) large enough so that~\((np)^{\frac{g}{g-c+1}} = n^{\frac{(1-\beta)g}{g-c+1}} \leq n^{1-\frac{\beta}{2}}\) and using~(\ref{e_cyc_est}) and~(\ref{e_deg_est}), we get that~\(\frac{\chi_{\eta,c}(G)}{n} \longrightarrow 0\) in probability. This completes the proof of Theorem~\ref{thm_main}\((b).\)~\(\qed\)

\subsection*{\em Acknowledgement}
I thank IMSc faculty for crucial comments and also thank IMSc for my fellowships.

\bibliographystyle{plain}

\end{document}